\documentclass[10pt]{article} 
\makeatletter
\usepackage{color,lineno}

\newtheorem {theorem}{Theorem}
\newtheorem {corol}{Corollary}
\newtheorem {consec}{Corollary}
\newtheorem {lemma}{Lemma}

\newcommand{\point}{\hspace{-1.75mm}{\bf.\ }}

\newcommand{\btheorem}{\begin{theorem}\point}
\newcommand{\etheorem}{\end{theorem}}

\newcommand{\blemma}{\begin{lemma}\point}
\newcommand{\elemma}{\end{lemma}}

\newcommand{\bpro}{\begin{pro}\point}
\newcommand{\epro}{\end{pro}}

\newcommand{\bcorol}{\begin{corol}\point}
\newcommand{\ecorol}{\end{corol}}

\newcommand{\bnote}{\begin{note}\point}
\newcommand{\enote}{\end{note}}

\newcommand{\bconsec}{\begin{consec}\point}
\newcommand{\econsec}{\end{consec}}

\newcommand{\bdefin}{\begin{defin}\point}
\newcommand{\edefin}{\end{defin}}
\newtheorem {pro}{Proposition}

\begin{document}
\title{On the existence of a $(2,3)$-spread in $V(7,2)$}
\author{Olof Heden and Papa A. Sissokho\thanks{Second author supported by grant KAW 2005.0098 
from the Knut and Alice Wallenberg Foundation.}} 

\maketitle
\begin{abstract}
An $(s,t)$-spread in a finite vector space $V=V(n,q)$ is a collection $\mathcal F$ of $t$-dimensional subspaces of $V$ with the property that every $s$-dimensional subspace of $V$ is contained in exactly one member of $\mathcal F$. It is remarkable that no $(s,t)$-spreads has been found yet, except in the case $s=1$.

In this note, the concept $\alpha$-point to a $(2,3)$-spread $\mathcal F$ in 
{$V=V(7,2)$} is introduced. A classical result of Thomas, applied to the vector 
space $V$, states that all points of $V$ cannot be $\alpha$-points to a given $(2,3)$-spread 
$\mathcal F$ in $V$. {In this note, we strengthened this result by proving that} 
every $6$-dimensional subspace of $V$ must contain at least one point that is not an $\alpha$-point to a given $(2,3)$-spread of $V$.
\end{abstract}

\section{Introduction}

An {\it $(s,t)$-spread} in the finite vector space $V=V(n,q)$ over GF$(q)$ is a collection 
$\mathcal F$ of $t$-dimensional subspaces of $V$ with the property that every $s$-dimensional subspace of $V$ is contained in exactly one member of $\mathcal F$. So far no $(s,t)$-spread, with $s>1$, has been found, and it was conjectured by Metsch that none {exists}, see \cite{metsch} for a survey. 

If there exists an $(s,t)$-spread $\mathcal F$ in $V$ then for any point $P$ in $V$, the members of $\mathcal F$ that contain $P$ induce an $(s-1,t-1)$-spread ${\mathcal F}_P$ in the quotient space $V/P$. A $(1,t)$-spread, or for short {\it spread}, ${\mathcal S}$ of $V$ is called {\em geometric} if for any three members $S_1$, $S_2$ and $S_3$ of $\mathcal S$ such that $S_3\cap\langle S_1\cup S_2\rangle\not=\{0\}$, we have  $S_3\subseteq\langle S_1\cup S_2\rangle$. 

Thomas~\cite{thomas} proved the following theorem.
\begin{theorem}\label{theorem:0}
Given a $(2,t)$-spread ${\mathcal F}$ of $V=V(n,q)$, there exists a point $P$ in 
$V$ such that the derived $(1,t-1)$-spread ${\mathcal F}_P$ is not geometric.
\end{theorem}

It must be remarked that geometric spreads are the spreads that are most natural and 
 ``easiest'' to construct, although most of the spreads are not geometric. 

The existence of $(2,3)$-spreads in $V(7,2)$ is the ``first'' open case for this conjecture. In this note, we give a property of $(2,3)$-spreads in $V(7,2)$, which, in this particular case, 
 yields the result of Thomas as a corollary.

Assume that $\mathcal F$ is a $(2,3)$-spread in $V=V(7,2)$. 
{As every spread in a $6$-dimensional subspace $U$ of $V$ is of size $21$, we get that every 1-dimensional subspace $P$, or {\it point}, of $V$ is contained in $21$ members of $\mathcal F$. As each of these $21$ members of $\mathcal F$ contains $7$ points, of which three belongs to $U$, it follows that $U$ contains 45 members of $\mathcal F$. Similarly, we may derive that every point $P$ in $U$ is contained in exactly $5$ of these 45 members of $\mathcal F$ and that every $5$-dimensional subspace $T$ of $U$ contains exactly five members of $\mathcal F$.} 

We will say that a point P is an {\it $\alpha$-point} to $\mathcal F$ if every $5$-dimensional subspace $T$ of $V$ that contains two of the members of $\mathcal F$ that meet at $P$, has the property that all its five members {from}  $\mathcal F$ will meet at the point $P$. From the definition of a geometric spread, it follows that in the case of $(2,3)$-spreads in $V=V(7,2)$,  Theorem~\ref{theorem:0} of Thomas states that at least one point of $V$ is not an $\alpha$-point to $\mathcal F$. 

We will show the following Theorem.
\begin{theorem}\label{theorem:2}
Assume that $\mathcal F$ is a $(2,3)$-spread in $V=V(7,2)$. Every $6$-dimensional subspace of $V$ contains at least one point which is not an $\alpha$-point to $\mathcal F$.
\end{theorem}

\section{Proof of Theorem \ref{theorem:2}}
Assume that $\mathcal F$ is a $(2,3)$-spread in $V=V(7,2)$. Let $U$ be any $6$-dimensional subspace of $V$. Assume that all points in $U$ are $\alpha$-points to $\mathcal F$. Then every $5$-dimensional subspace $T$ of $U$ will contain a point $P$ where all its five members of $\mathcal F$ meet. This point $P$ will be called the {\it $\alpha$-point of $T$}. Moreover, each point $P$ of $U$ is contained in exactly five of the members of $\mathcal F$ that belong to $U$, and hence these five members of $\mathcal F$ that meet the point $P$ will all belong to the same $5$-dimensional subspace $T$ of $U$.

{{We claim that} there is a {$4$-dimensional} subspace $W$ of $U$ that does not contain any member of $\mathcal F$. To see this, just observe that every $3$-dimensional subspace of a $5$-dimensional subspace $T$ of $U$ is contained in exactly three $4$-dimensional subspaces of $T$, and as $T$ contains exactly five members of $\mathcal F$, there will be at least 16 subspaces $W$ of dimension $4$ of $T$ that do not contain any member of 
$\mathcal F$.}
Such a $4$-dimensional subspace $W$ of $U$ will be called a {\it poor space}. 

There are three $5$-dimensional subspaces $T_1$, $T_2$ and $T_3$ of $U$ such that
\begin{equation}\label{eq:2}
W=T_1\cap T_2=T_1\cap T_3=T_2\cap T_3\;,\qquad\hbox{and}\qquad U=T_1\cup T_2\cup T_3\;.
\end{equation} 
For $1\leq i\leq3$, let $P_i$ be the $\alpha$-point in the space $T_i$. 

We first note that none of the points $P_1$, $P_2$, or $P_3$ belongs to $W$. 

To prove this fact, assume for instance that $P_1$ belongs to $W$.  Since $W$ is a poor $4$-dimensional space, each of the five members of $\mathcal F$ that belongs to $U$ and contains the point $P_1$ meet $W$ in two points, besides the point $P_1$.
This leads to a contradiction since $W$ contains $15$ points and every point $Q\neq P_1$ in $T_1$
(and thus in $W$) belongs to exactly one of the five members of $\mathcal F$ in $U$ that meet the point $P_1$.

Since $\mathcal F$ is a $(2,3)$-spread and since the points $P_i$, $1\leq i\leq3$, do not belong to $W$ and they are the $\alpha$-points of the respective spaces $T_i$, we can conclude that the members of $\mathcal F$ that are subspaces of $T_i$ will intersect $W$ in a spread ${\mathcal S}_i$. Furthermore, since $\mathcal F$ is a $(2,3)$-spread, these three spreads are mutually disjoint.

Now, let $Q$ be any point of $W$. Let $T_Q$ denote the unique $5$-dimensional subspace of $U$, that contains the two members of $\mathcal F$ that meet the point $Q$ and belong to $T_1$ and $T_2$, respectively. We note from Equation (\ref{eq:2}) that $P_1\not\in T_2\cup T_3$ and $P_2\not\in T_1\cup T_3$. {Hence, $T_Q$ cannot be one of the spaces $T_i$, $1\leq i\leq3$. As these are the only $5$-dimensional subspaces of $U$ that contain $W$, it follows that
\[
\dim(T_Q\cap W)\leq3\;.
\]   }

Moreover, since all $5$-dimensional subspaces of $U$ have a {unique} point where all its members of $\mathcal F$ meet, and as there are two members of $\mathcal F$ in $T_Q$ meeting 
$Q$, we conclude that $Q$ is the $\alpha$-point of the space $T_Q$. This implies that the member of $\mathcal F$ that is a subspace of $T_3$ and meets the point $Q$ must also belong to $T_Q$. This space will be denoted by $Z_{Q,3}$; and we define $Z_{Q,1}$ and $Z_{Q,2}$ similarly. For 
$1\leq i\leq3$, the intersection of $Z_{Q,i}$ with $W$ is a $2$-dimensional subspace which we denote by 
$L_{Q,i}$. 

Now, the space $Z_{Q,3}$ is completely contained in $T_Q$ and intersects $W$ in the $2$-dimensional space $L_{Q,3}$, which thus also must be a subspace of $T_Q$, so,
\begin{equation}\label{eq:1}
L_{Q,3}\subseteq T_Q\cap W=\langle L_{Q,1},L_{Q,2}\rangle\;.
\end{equation}
The last step in our proof is to show that there is at least one point $Q$ in $W$, for which the above relation does not hold.

{Let us assume for a moment that}
\[
{\mathcal S}_1=\{\;L_1, L_2, \dots, L_5\;\}\qquad\hbox{and}\qquad{\mathcal S}_2=\{\;L_1', L_2', \dots, L_5'\;\}\;.
\]
Every member{, or {\it line},} of ${\mathcal S}_2$ intersects three members of ${\mathcal S}_1$. 
{Without loss of generality, we may assume that the line $L_5'$ does not intersect the lines} {$L_1$ and $L_2$. These two lines together contain $6$ points. Each of these $6$ points is contained in {exactly one} of the lines of ${\mathcal S}_2$. As a line contains $3$ points we get that there must be two lines, say $L_1'$ and $L_2'$, of ${\mathcal S}_2$ that meet both $L_1$ and $L_2$.

Let $Q=L_1\cap L_1'$, $Q'=L_2\cap L_2'$, $R_1=L_1\cap L_2'$ and $R_2=L_2\cap L_1'$, i.e., with the original notation }
\begin{equation}\label{eq:3}
L_{Q,1}\cap L_{Q',2}=R_1\qquad\hbox{and}\qquad L_{Q,2}\cap L_{Q',1}=R_2\;.
\end{equation}
{Then the line $L$, that meets the points $R_1$ and $R_2$, satisfies the following relation}  
\[
L=\langle R_1,R_2\rangle=(T_Q\cap W)\cap(T_{Q'}\cap W)\;.
\]
If the relation (\ref{eq:1}) holds for all points $Q$ of $W$, then $L$ will meet both the spaces $L_{Q,3}$ and $L_{Q',3}$. Note that $L$ contains just three points, the above defined two points $R_1$ and $R_2$, and a third point $R_3$. So from Equation (\ref{eq:3}), we can infer that both the spaces $L_{Q,3}$ and $L_{Q',3}$ must meet $L$ at the point $R_3$. This contradicts the fact that ${\mathcal S}_3$ is a spread and the proof is complete.

\bigskip

\noindent
O. Heden, Department of Mathematics, KTH, S-100 44 Stockholm, Sweden.

\noindent
(olohed@math.kth.se)

\bigskip
\noindent
P. A. Sissokho, 4520 Mathematics Department, Illinois State University, 
Normal, Illinois 61790--4520, U.S.A. 

\noindent
(psissok@ilstu.edu)
\end{document}